\documentclass{amsart}
\usepackage{hyperref}

\usepackage{amssymb} 
\usepackage{enumitem}
\usepackage{multirow}
\setlist[itemize]{leftmargin=12mm}
\setlist[enumerate]{leftmargin=12mm}
\usepackage{comment}
\usepackage[all,cmtip]{xy}
\usepackage{url}

\DeclareMathOperator{\End}{End}

\DeclareMathOperator{\GL}{GL}

\DeclareMathOperator{\SL}{SL}

\DeclareMathOperator{\Gal}{Gal}

\newcommand{\Q}{{\mathbb Q}}
\newcommand{\Z}{{\mathbb Z}}
\newcommand{\F}{{\mathbb F}}

\newcommand{\galas}{\Gal(\Qb/\Q)}

\newcommand{\abc}[4]{\left(\begin{smallmatrix}#1&#2\\#3&#4\end{smallmatrix}\right)}

\def\C{{\ensuremath{\mathbb{C}}}}

\def\Qb{{\ensuremath{\overline{\mathbb{Q}}\,}}}
\def\O{{\ensuremath{\mathcal{O} }}}

\def\F{{\ensuremath{\mathbb{F}}}}

\def\p{{\ensuremath{\mathfrak{p}\,}}}
\def\T{{\ensuremath{\mathbb{T}}}}

\def\nc{\mathrm{N}}

\DeclareMathOperator{\Tr}{Trace}

\DeclareMathOperator{\Frob}{Frob}
\DeclareMathOperator{\cusps}{Cusps}
\DeclareMathOperator{\tate}{Tate}
\DeclareMathOperator{\can}{can}

\DeclareMathOperator{\cohom}{H}

\DeclareMathOperator{\res}{Res}
\DeclareMathOperator{\spec}{Spec}

\begin{document}

\newtheorem{thm}{Theorem}[section]
\newtheorem{lem}{Lemma}[section]
\newtheorem{prop}[lem]{Proposition}
\newtheorem{cor}[thm]{Corollary}
\newtheorem{conj}{Conjecture}

\theoremstyle{definition}
\newtheorem*{remark}{Remark}
\newtheorem*{ex}{Example}

\theoremstyle{remark}
\newtheorem{ack}{Acknowledgement}

\title[The minimal level of realization for a mod~$\ell$ eigenvalue system]{A note on the minimal level of realization for a mod~$\ell$ eigenvalue system}
\author{Samuele Anni}

\address{Mathematics Institute\\
	University of Warwick\\
	Coventry\\
	CV4 7AL \\
	United Kingdom}
\email{samuele.anni@gmail.com}

\date{\today}
\thanks{The author is supported by EPSRC Programme Grant 
\lq LMF: L-Functions and Modular Forms\rq\  EP/K034383/1.
}


\begin{abstract}
In this article we give a criterion for a mod~$\ell$ eigenvalue system attached to a  mod~$\ell$ Katz cuspform to arise from lower level or weight. 
Namely, we prove the following: the eigenvalue system associated to a ring homomorphism  $f:\T \to \overline{\F}_\ell$ 
from the Hecke algebra of level $\Gamma_1(n)$ and weight $k$ to $\overline{\F}_\ell$, where $\ell$ is a prime not dividing $n$ and  $1\leq k \leq \ell +1$,
arises from lower level or weight if there exists a prime $r$ dividing $n\ell$ such that  
$$ \dim_{\overline{\F}_\ell} \bigcap_{p \neq r} \ker \left( T_p-f(T_p), S(n,k)_{\overline{\F}_\ell}\right)>1,$$
where $T_p$ is the $p$-th Hecke operator and $S(n,k)_{\overline{\F}_\ell}$ is the space of  mod~$\ell$ Katz cuspforms of level $\Gamma_1(n)$ and weight $k$. 
\end{abstract}
\maketitle

\section{Introduction and preliminaries}
In this article we analyse degeneracy maps between modular curves in positive characteristic in order to deduce statements on residual modular Galois representations and their level of realization. 
In \cite{edix1}, the study of the Hasse invariant and the Frobenius map lead to similar result about the minimal weight at which a mod~$\ell$ eigenvalue system attached to a Katz cuspform arises.

This article is based on the theory of Katz modular forms (\cite{k1}, \cite{k2}, \cite{edix3}) and modular curves (\cite{katzmaz}, \cite{gross}). 
Analogous results in characteristic zero of the some presented in this article are well-known (\cite[Section~5.7]{dish}). Anyway, in positive characteristic they can be obtained through the study of modular curves over finite fields and a geometric interpretation of degeneracy maps. 
The statements that we will prove do not rely on lifting to characteristic zero, and so, in particular, they hold also in weight $1$. 
The results herein may lead to a formulation of multiplicity one statements in positive characteristic, anyway, this will be subject of a subsequent paper.

For any pair of positive integers $n$ and $k$, called respectively level and weight, we define the following objects: $M(n,k)_\C$, 
the complex vector space of modular forms of weight $k$ on $\Gamma_1(n)$ and $S(n,k)_\C$, the subspace of cuspforms. 
Let $\T$ be the associated Hecke algebra, i.e.\ the $\Z$-subalgebra of $\End_{\C}(S(n,k)_\C)$ generated by the Hecke operators $T_p$ for every prime $p$ and the diamond operators $\left\langle d\right\rangle$ for every $d \in (\Z/n\Z)^{\ast}$. 
The Hecke algebra $\T$ is finitely generated as a $\Z$-module (\cite[Theorem~2.5.11]{edix2}, \cite[p.234]{dish}). 

Let us fix an algebraic closure $\overline{\Q}$ of $\Q$, we will denote as $G_\Q$ the absolute Galois group $\Gal(\overline{\Q}/\Q)$.
To any ring homomorphism from the Hecke algebra $\T$ of level $\Gamma_1(n)$ and weight $k$ to a finite field of characteristic not dividing $n$ 
we associate a continuous semisimple Galois representation thanks to the following result due to Shimura and Deligne (\cite[Th\'eor\`eme~6.7]{ds}, \cite[Theorem~2.5.2]{edix2}):

\begin{thm}[Deligne, Shimura]\label{maxrep}
Let $n$ and $k$ be positive integers. Let $\F$ be a finite field of characteristic $\ell$ and $f:\T \to \F$ a morphism of rings. Then there is a continuous semi-simple representation $\rho_f: G_\Q \to \GL_2(\F)$ that is unramified outside $n\ell$ such that for all primes $p$ not dividing $n\ell$ we have:
$$\Tr(\rho_f(\Frob_p)) = f(T_p)\mbox{ and }\det(\rho_f(\Frob_p)) = f(\left\langle p\right\rangle)p^{k-1}\mbox{ in } \F.$$
Such a $\rho_f$ is unique up to isomorphism.
\end{thm}
A ring morphism as in Theorem~\ref{maxrep} corresponds to an eigenform with coef\-ficients in $\F$, for more details see \cite[p.58]{edix2}. 
We will denote by $\nc(\rho_f)$ the Artin conductor away from $\ell$ of $\rho_f$ and by $k(\rho_f)$ the weight of $\rho_f$, see \cite{edix3}.

The space $S(n,k)_\C$ is decomposed by characters, let $S(n, k, \tilde{\epsilon})_\C$ denote the space of cuspforms with character $\tilde{\epsilon}$. 
We associate a Hecke algebra to each $S(n, k, \tilde{\epsilon})_\C$, as before: this Hecke algebra, that we will denote as $\T_{\tilde{\epsilon}}$, 
is the $\O_{\tilde{\epsilon}}$-subalgebra in the complex endomorphism ring of $S(n,k,\tilde{\epsilon})_\C$ generated by the Hecke operators $T_n$ for $n\geq 1$, 
where $\O_{\tilde{\epsilon}}$ is ring of integers of the number field containing the image of $\tilde{\epsilon}$. 
Note that the algebra $\T_{\tilde{\epsilon}}$ is free of finite rank as a $\Z$-module. From Theorem~\ref{maxrep} it follows that:
\begin{cor}  
Let $n$ and $k$ be positive integers, let $\tilde{\epsilon}:(\Z/n\Z)^\ast \to \C^\ast$ be a character and let 
$f : \T_{\tilde{\epsilon}} \to \overline{\F}_\ell$ be a morphism of rings. Then there is a continuous semi-simple representation 
$\rho_f : G_\Q \to \GL_2(\overline{\F}_\ell)$ that is unramified outside $n\ell$, such that for all primes $p$ not dividing $n\ell$ we have:
$$\Tr(\rho_f(\Frob_p)) = f(T_p)\mbox{ and }\det(\rho_f(\Frob_p)) = f(\left\langle p\right\rangle)p^{k-1}\mbox{ in } \overline{\F}_\ell.$$
Such a $\rho_f$ is unique up to isomorphism.
\label{obj}
\end{cor}
A ring morphism as above corresponds to a mod~$\ell$ eigenform with character $\epsilon:(\Z/n\Z)^\ast \to \overline{\F}_\ell$ given by  
$a \mapsto \epsilon(a):= f(\left\langle a\right\rangle)$ for $a\in (\Z/n\Z)^{\ast}$. Equivalently, $\det(\rho_f)=\epsilon\cdot \chi_\ell^{k-1}$ where 
$\chi_\ell$ is the cyclotomic character mod~$\ell$.
We will abuse of notation and denote by $f$ the mod~$\ell$ cuspidal eigenform corresponding to the ring homomorphism $f$ as in the above corollary.

Let us recall that the Hecke algebras $\T$ and $\T_{\tilde{\epsilon}}$ for level $n$, weight $k\geq 2$ and character 
$\tilde{\epsilon}:(\Z/n\Z)^\ast\to\O_{\tilde{\epsilon}}$, can be computed in polynomial time in $n$ and $k$, using deterministic algorithms 
based on mo\-dular symbols (\cite[Theorem~9.23 and Theorem~9.22]{ste}.

Let $\ell$ be a prime not dividing $n$ and $\overline{\F}_\ell$ an algebraic closure of $\F_\ell$, 
let $M(n,k)_{\overline{\F}_\ell}$ be the space of mod~$\ell$ Katz modular forms of weight $k$ for $\Gamma_1(n)$. 
Let $S(n,k)_{\overline{\F}_\ell}$ be the cuspidal subspace of $M(n,k)_{\overline{\F}_\ell}$, i.e\ the subspace whose elements have $q$-expansions that are power series with constant term zero at all cusps. 
For every character $\epsilon:(\Z/n\Z)^\ast \to \overline{\F}_\ell^\ast$, let $S(n,k,\epsilon)_{\overline{\F}_\ell}$ be the space of Katz cuspforms of weight $k$ for $\Gamma_1(n)$ and character $\epsilon$.
If ${n{>}4}$ then $S(n,k)_{\overline{\F}_\ell}$ is isomorphic to $\cohom^0(X_1(n)_{\overline{\F}_\ell},\omega^{\otimes k}(-\cusps))$, the space of global sections of the line bundle 
$\omega^{\otimes k}(-\cusps)$, while the space of Katz modular forms $M(n,k)_{\overline{\F}_\ell}$ is isomorphic to 
$\cohom^0(X_1(n)_{\overline{\F}_\ell},\omega^{\otimes k})$, see \cite[Proposition~$2.2$]{gross}, where $X_1(n)_{\overline{\F}_\ell}$ is the base change of the modular curve $X_1(n)$ over $\overline{\F}_\ell$.

The main result of this article is the following:
\begin{thm}\label{rid}
Let $n$ and $k$ be positive integers. Let $\ell$ be a prime not dividing $n$, such that $1\leq k \leq \ell +1$. 
Let $f:\T_{\tilde{\epsilon}} \to \overline{\F}_\ell$ be a ring homomorphism from the Hecke algebra of level $\Gamma_1(n)$, weight $k$ and character 
$\tilde{\epsilon}: \left(\Z/n\Z\right)^\ast \to \C^\ast$, to $\overline{\F}_\ell$. Let $T_p \in \T_{\tilde{\epsilon}}$ be the $p$-th Hecke operator and 
$\epsilon: \left(\Z/n\Z\right)^\ast \to \overline{\F}_\ell^\ast$ be the character defined by 
$a \mapsto f(\left\langle a\right\rangle)$ for all $a\in \left(\Z/n\Z\right)^\ast$. 

Assume that there exists a prime $r$ dividing $n\ell$ such that the $\overline{\F}_\ell$-vector space 
$$ V\colon = \bigcap_{p \neq r} \ker \left( T_p-f(T_p), S(n,k, \epsilon)_{\overline{\F}_\ell}\right),$$
has dimension bigger than $1$. 

If $r$ is different from $\ell$ then there exists $g \in S(n/r, k, \res(\epsilon))_{\overline{\F}_\ell}$ such that $f$ belongs to the span of 
$B^\ast_r g$ and $\alpha^\ast g$, where $B^\ast_r$ and $\alpha^\ast$ are the maps induced by the degeneracy maps.

If $r=\ell$ then $k\in\left\{\ell,\ell{+}1 \right\}$. 
In the first case there exists $g\in S(n, 1, \epsilon)_{\overline{\F}_\ell}$ such that $F(g)=f$, where $F$ is the Frobenius morphism. 
In the second case there exists $g\in S(n, 2, \epsilon)_{\overline{\F}_\ell}$ such that $A_\ell g=f$, where $A_\ell$ is the Hasse invariant.
\end{thm}

From the previous theorem we deduce the following results:
\begin{cor} 
The Galois representation $\rho_f: \galas \to \GL_2(\overline{\F}_\ell)$ attached to $f$, under the same hypotheses of Theorem~\ref{rid}, arises from a lower level or weight.
\end{cor}
\begin{cor} 
Under the same notation of Theorem~\ref{rid}, if $V$ has dimension $1$ then $f$ is of minimal weight and level, i.e. $\nc(\rho_f)=n$ and $k(\rho_f)=k$.
\end{cor}
We say that a form is minimal with respect to weight and level if the associate mod~$\ell$ eigenvalue system does not arise from forms of lower level or weight.
\section{Notation}

Let $n$ and $k$ be positive integers, let $\ell$ be a prime not dividing $n$ and $\epsilon$ a character over $\overline{\F}_\ell$ with conductor dividing $n$. 

Let us assume that $n=m p^r$ with $r\geq 1$ and $p$ not dividing $m$. 
We split the level structure according to this factorization, and denote the space of mod~$\ell$ cuspforms as 
$S(m,p^r,k,\overline{\epsilon})_{\overline{\F}_\ell}$. We will use a\-na\-lo\-gous notation for modular curves. 

On the modular curve $X_1(n)_{\overline{\F}_\ell}$ there are two degeneracy maps at $p$: we will denote them as $B_p$ and $\alpha$:
\begin{equation}\label{dege}
\xymatrix{
\quad & X_1(m, p^r)_{\overline{\F}_\ell} \ar[dl]_{B_p} \ar[dr]^\alpha & \quad \\
X_1(m, p^{r-1})_{\overline{\F}_\ell} & \quad & X_1(m, p^{r-1})_{\overline{\F}_\ell}.}
\end{equation}
Let us recall the description of the degeneracy maps using the moduli interpretation for $X_1(n)_{\overline{\F}_\ell}$. 
Let $E/S$ be an el\-liptic curve over an $\overline{\F}_\ell$-scheme $S$, with $P$ and $Q$ respectively points of order $m$ and $p^r$. 
The map $\alpha$ is the forgetful map $\alpha\colon(E,P,Q)\longmapsto (E, P, pQ)$ and the map $B_p$ is the $p$-th degeneracy map defined by 
$B_p \colon (E,P,Q)\mapsto (E/\langle p^{r-1} Q\rangle, \beta(P), \beta(Q))$, 
where 
$$\langle p^{r-1} Q\rangle\, \rightarrowtail \;E \stackrel{\beta}{\longrightarrow} E/\langle p^{r-1} Q\rangle$$
is the degree $p$ isogeny whose kernel is generated by $p^{r-1} Q$.

Let us fix $\zeta_{p^r}$ a root of unity of order $p^r$ in $\overline{\F}_\ell$. 
The Atkin-Lehner involution $w_{\zeta_{p^r}}$ on $X_1(n)_{\overline{\F}_\ell}$ is defined as follows:
\begin{equation}\label{atkleh}
\begin{array}{lcl}
X_1(m, p^r)_{\overline{\F}_\ell}  & \stackrel{w_{\zeta_{p^r}}}{\longrightarrow} &  X_1(m, p^r)_{\overline{\F}_\ell}\\
(E,P,Q) & \longmapsto & (E/\langle Q\rangle, \beta_r(P), Q')
\end{array}
\end{equation}
where $\beta_r$ is the degree $p^r$ isogeny whose kernel is generated by $Q$, and $Q'$ belongs to $\ker(\beta_r^t)(S)$:
$$\xymatrix{\langle Q\rangle\, \ar@{>->}[r] & E \ar@/^/[rr]|{\beta_r} && E/\langle Q\rangle \ar@/^/[ll]|{\beta_r^{t}} & \,\langle Q'\rangle\ar@{>->}[l]}$$
and it is the unique element such that $e_\beta(Q,Q')=\zeta_{p^r}$, where $e_\beta$ is the Weil pairing related to $\beta$ i.e.\ the perfect $\mu_{p^r}$ valued pairing between $\ker(\beta)$ and $\ker(\beta^t)$, see \cite[Section~2.8]{katzmaz}.

Let us recall basic facts about Tate curves, see \cite[Chapter~VII]{derap} and \cite[Chapter~V]{silv2}. 
For every positive integer $d$, the $d$-th Tate curve is a certain generalized elliptic curve $\tate(q^d) \to \spec\Z[[q]]$ that becomes a N\'eron $d$-gon after base change to the zero locus of $q$ and that is a smooth elliptic curve over $\spec\Z[[q]][q^{-1}]$, the complement of this zero locus.
Let $d$, $e$ and $n$ be positive integers such that $n$ is the least common multiple of $d$ and $e$. 
Then the curve $\tate(q^d)$  over $\spec\Z[[q, \zeta_e]]$, where $\zeta_{e}$ is a $e$-th root of unity, admits a $\Gamma_1(n)$-structure. Each choice of $d$, $e$ and of the a $\Gamma_1(n)$-structure gives rise to an injective map from $M(\Gamma_1(n),k)_{\overline{\F}_\ell}$ to $\overline{\F}_\ell \otimes_\Z \Z[[q, \zeta_e]]$, which is called the $q$-expansion map relative to $\tate(q^d)$ with the given $\Gamma_1(n)$-structure.
Evaluating any mod~$\ell$ modular form in $M(\Gamma_1(n),k)_{\overline{\F}_\ell}$ on the Tate curve $\tate(q^n)$ over $\spec\Z[[q]]$ with respect to the 
$\Gamma_1(n)$-structure given by $\psi(i) = q^i$ for $i \in \Z/n\Z$, we obtain a particular $q$-expansion called the $q$-expansion at $\infty$ (\cite[Sections~8.8 and 8.11]{katzmaz}). 
For any $f \in M(\Gamma_1(n),k)_{\overline{\F}_\ell}$ we define $a_m(f)$ to be the $m$-th coefficient of this $q$-expansion.

\section{Level lowering for Katz cuspforms}
\begin{lem} \label{old1}
Let $n$ and $k$ be positive integers and let $p$ be a prime strictly dividing $n$ with $n/p>4$. Let $\ell$ be a prime not dividing $n$ and such that $1\leq k \leq \ell{+}1$. 
Let $f\in S(n, k)_{\overline{\F}_\ell}$ such that 
\begin{equation*}
f(\tate(q), \zeta_{n/p}, \zeta_{p})= \sum_{j\geq 1} a_j(f) q^{j} \left(\frac{dt}{t}\right)^{\otimes k}\in \overline{\F}_\ell\left[\left[q^p\right]\right]\left(\frac{dt}{t}\right)^{\otimes k},
\end{equation*}
where $\zeta_{n/p}$ and $\zeta_p$ are respectively $n/p$-th and $p$-th roots of unity. 
Then there exists a unique cuspform $g\in S(n/p,k)_{\overline{\F}_\ell}$, such that $B_p^\ast g=f$, where $B_p^\ast$ is the pullback of 
$B_p\colon X_1(n)_{\overline{\F}_\ell} \to X_1(n/p)_{\overline{\F}_\ell}$, the $p$-th degeneracy map.
\end{lem}
\proof Since $p$ strictly divides $n$, let us write $n=mp$. By hypothesis, the $q$-expansion of $f$ at the cusp $\infty$ is such that $a_j(f)=0$ for all positive integers $j$ not divisible by $p$. 
As $m>4$ by hypothesis, $f$ is an element of $\cohom^0(X_1(n)_{\overline{\F}_\ell}, \omega^{\otimes k}(-\cusps))$. From the description of the degeneracy maps given in the previous section, it is easy to show that the 
following diagram is commutative: 
$$\xymatrix{
X_1(n)_{\overline{\F}_\ell}  \ar[d]_{B_p} \ar[r]^{w_{\zeta_p}} &  X_1(n)_{\overline{\F}_\ell}  \ar[d]^\alpha\\
X_1(m)_{\overline{\F}_\ell} \ar@{=}[r] & X_1(m)_{\overline{\F}_\ell} }$$
where $B_p$ and $\alpha$ are the degeneracy maps defined in (\ref{dege}) and $w_{\zeta_p}$ is the Atkin-Lehner map associated to the isogeny $\beta$, defined in (\ref{atkleh}). The isogeny $\beta$ induces a morphism $\beta^\ast$ from the line bundle $\omega^{\otimes k}(-\cusps)$ on $X_1(m)_{\overline{\F}_\ell}$ to the line bundle $\omega^{\otimes k}(-\cusps)$ on $X_1(n)_{\overline{\F}_\ell}$, since for each 
$(E/S, P, Q)$ we have an isomorphism of invertible $O_S$-modules $\beta^\ast$ between $\omega_{E/\langle Q\rangle}$ and $\omega_E$. The pullback $w^\ast_{\zeta_p}: S(n, k)_{\overline{\F}_\ell}\to S(n, k)_{\overline{\F}_\ell}$ is defined by 
$$(w^\ast_{\zeta_p} h)(E, P, Q)= \beta^\ast( h(w_{\zeta_p} (E, P, Q)))=\beta^\ast( h(E/\langle Q\rangle, \beta(P), Q')),$$ 
where $h$ belongs to $S(n, k)_{\overline{\F}_\ell}$, see \cite[Section~6]{gross}. 
Therefore, there exists a unique $f' \in  S(n, k)_{\overline{\F}_\ell}$ satisfying $w^\ast_{\zeta_p} f' =f$. 
Indeed, the $q$-expansion at the cusp $\infty$ gives:
\begin{eqnarray*}
(w^\ast_{\zeta_p} f')(\infty) &=& (w^\ast_{\zeta_p} f')(\tate(q), \zeta_{m}, \zeta_{p})= \beta^\ast( f'(\infty'))=  \\
&=& \beta^\ast( f'(\tate(q^p), \zeta_{m}^p, q))=\beta^\ast\left(\sum_{j\geq 1} a_j(f')q^j \left(\frac{dt}{t}\right)^{\otimes k}\right)=\\
&=& \left(\sum_{j\geq 1} a_j(f')q^j \right)p^k\left(\frac{dt}{t}\right)^{\otimes k},
\end{eqnarray*} where  $\infty':= w_{\zeta_p} (\infty)$ and the last equality follows from:
$$\xymatrix{\mu_p \; \ar@{>->}[r] & \tate(q) \ar@{>>}[r]^\beta & \tate(q^p)},$$
hence $\beta^\ast(dt/t)=p (dt/t)$. Since 
$f(\tate(q), \zeta_{m}, \zeta_{p})\in \overline{\F}_\ell\left[\left[q^p\right]\right]\left(dt/t\right)^{\otimes k}$, we deduce that for all positive integer $j$ we have $a_j(f')p^k=a_j(f)$ and  
$$f'(\tate(q^p), \zeta_{m}^p, q)\in \overline{\F}_\ell\left[\left[q^p\right]\right]\left(\frac{dt}{t}\right)^{\otimes k}.$$
In order to prove that $f$ comes from $X_1(m)_{\overline{\F}_\ell}$ via $B_p^\ast$, it is enough to show that $f'$ does via $\alpha$. 
To prove this claim, consider the following covering of modular curves: 
$$\xymatrix{
X(\Gamma_1(m),\Gamma(p)^{\zeta_p - \can})_{\overline{\F}_\ell}  \ar[d]_{\gamma} & (E,P,Q_1,Q_2)  \ar@{|->}[d] & 
(\tate(q^p),\zeta_{m}^p, q, \zeta_p) \ar@{|->}[d]\\
X(\Gamma_1(m),\Gamma_1(p))_{\overline{\F}_\ell}  \ar[d]_{\alpha} & (E,P,Q_1)  \ar@{|->}[d] & (\tate(q^p),\zeta_{m}^p, q) \ar@{|->}[d]\\
X(\Gamma_1(m))_{\overline{\F}_\ell}   & (E,P)  & (\tate(q^p),\zeta_{m}^p)}$$
where $Q_1, Q_2$ are such that $e_p(Q_1, Q_2)=\zeta_p$, and $\alpha, \gamma$ are both forgetful maps. 
On $X(\Gamma_1(m),\Gamma(p)^{\zeta_p - \can})_{\overline{\F}_\ell} \to X(\Gamma_1(m))_{\overline{\F}_\ell}$ there is the natural $\SL_2(\F_p)$ action, given by the action of $\SL_2(\F_p)$ on $\Gamma(p)^{\zeta_p - \can}$. 
This action is trivial on $X(\Gamma_1(m))_{\overline{\F}_\ell}$: 
$$\xymatrix{
X(\Gamma_1(m),\Gamma (p)^{\zeta_p - \can})_{\overline{\F}_\ell}  \ar[d]^{\left({\begin{smallmatrix}1 & \ast  \\0 & 1\end{smallmatrix}}\right)} 
\ar@/_6pc/ [dd]_{\SL_2(\F_p)} \\
X(\Gamma_1(m),\Gamma_1(p))_{\overline{\F}_\ell}  \ar[d] \\
X(\Gamma_1(m))_{\overline{\F}_\ell}}$$
Therefore, there exists $g \in  S(m,k)_{\overline{\F}_\ell}$ such that $B^\ast_p g=f$ if and only if 
$\gamma^\ast f'$ is $\SL_2(\F_p)$-invariant.

The map $\gamma : X(\Gamma_1(m),\Gamma(p)^{\zeta_p - \can})_{\overline{\F}_\ell} \to X(\Gamma_1(m),\Gamma(p))_{\overline{\F}_\ell}$ corresponds to the action of elements of the form 
$\left({\begin{smallmatrix}1 & \ast  \\0 & 1\end{smallmatrix}}\right)\in \SL_2(\F_p)$ and $\gamma^\ast f'$ is invariant under the action of 
$\left({\begin{smallmatrix}1 & \ast  \\0 & 1\end{smallmatrix}}\right)$ by construction. Since $\left\{\left({\begin{smallmatrix}1 & 0 \\\ast & 1\end{smallmatrix}}\right),\left({\begin{smallmatrix}1 & \ast  \\0 & 1\end{smallmatrix}}\right)\right\}$ generate $\SL_2(\F_p)$, and $\left({\begin{smallmatrix}1 & 0 \\1 & 1\end{smallmatrix}}\right)$ generates the subgroup $\left({\begin{smallmatrix}1 & 0 \\\ast & 1\end{smallmatrix}}\right)$ in $\SL_2(\F_p)$, in order to conclude it is enough to prove that $\gamma^\ast f'$ is invariant under the action of $\left({\begin{smallmatrix}1 & 0 \\1 & 1\end{smallmatrix}}\right)$. 

Let us remark that 
$\left({\begin{smallmatrix}1 & 0 \\1 & 1\end{smallmatrix}}\right)$ fixes 
$\infty'\in X(\Gamma_1(m),\Gamma (p)^{\zeta_p - \can})_{\overline{\F}_\ell}$ and it acts on $\overline{\F}_\ell\left[\left[q\right]\right]$ as
$$
\left\{
\begin{array}{l}
\mbox{ identity on }\overline{\F}_\ell\\
q \mapsto \zeta_p q	
\end{array}
\right..
$$  
It follows that $\left({\begin{smallmatrix}1 & 0 \\1 & 1\end{smallmatrix}}\right)$ acts on 
$\overline{\F}_\ell\left[\left[q\right]\right]\left(dt/t\right)^{\otimes k}$ mapping $t \mapsto t$, by definition, and so 
$\left(dt/t\right)^{\otimes k}\mapsto\left(dt/t\right)^{\otimes k}$. 
As $a_j(f')=0$ for all po\-si\-ti\-ve integers $j$ not divisible by $p$, then $\left({\begin{smallmatrix}1 & 0 \\1 & 1\end{smallmatrix}}\right)^\ast f'$ has the same $q$-expansion as $f'$. 
Since  
$X(\Gamma_1(m),\Gamma (p)^{\zeta_p - \can})_{\overline{\F}_\ell}$ is an irreducible and reduced curve, this means that 
$\left({\begin{smallmatrix}1 & 0 \\1 & 1\end{smallmatrix}}\right)^\ast\gamma^\ast  f'=\gamma^\ast f'$.

The cuspform $g \in S(m,k)_{\overline{\F}_\ell}$ such that $B^\ast_p g=f$ is unique: from the one hand, pullback does not change $q$-expansions and, on the other hand, the $q$-expansion principle holds (\cite[Section~1.12]{k1}). 
\endproof

Now we want to generalize the previous lemma for levels $n$ and primes $p$ such that $n$ is highly divisible by $p$. In order to do so, we need some preliminary results which will help us understanding the modular curves covering and the interaction between level structure and Atkin\--Lehner operators.
\begin{lem}\label{gen}
Given $r \in \Z_{>0}$ and a prime $p$, the set 
$$\left\{\left({\begin{matrix}1 & 0 \\p^r & 1\end{matrix}}\right),\left({\begin{matrix}1 & 1  \\0 & 1\end{matrix}}\right)\right\}\subset \SL_2(\Z/p^{r+1}\Z)$$ 
generates the subgroup 
$$\overline{\Gamma_1(p^r)}\colon= 
\left\{\left(\begin{matrix}1+p^r a  & b\\p^r  c & 1+p^r  d \end{matrix}\right)\in \SL_2(\Z/p^{r+1}\Z)\right\}.$$ 
\end{lem}
\proof This is a straightforward computation: every matrix in $\overline{\Gamma_1(p^r)}$ can be written as a finite product of 
$\left({\begin{smallmatrix}1 & 0 \\p^r & 1\end{smallmatrix}}\right)$ and $\left({\begin{smallmatrix}1 & 1  \\0 & 1\end{smallmatrix}}\right)$. 
Indeed, let $\left({\begin{smallmatrix}1+p^r a & b\\p^r c & 1+p^r d\end{smallmatrix} }\right)$ in $\SL_2(\Z/p^{r+1}\Z)$, with $a,b,c,d \in \F_p$, 
then 
$$ \left({\begin{matrix}1+p^r a & b\\p^r c & 1+p^r d\end{matrix} }\right)= 
\left({\begin{matrix}1 & 0 \\p^r & 1\end{matrix}}\right)^{c-1} \left({\begin{matrix}1 & 1  \\0 & 1\end{matrix}}\right)^a
\left({\begin{matrix}1 & 0 \\p^r & 1\end{matrix}}\right)^{-1}
\left({\begin{matrix}1 & 1  \\0 & 1\end{matrix}}\right)^{(p^r a-1)(a-b)}.$$
\endproof

\begin{lem}\label{isog}
Let $E/S$ be an elliptic curve over an $\overline{\F}_\ell$-scheme $S$, let $p$ be a prime different from $\ell$. Let $r\geq 1$ be an integer and  
$$\xymatrix{ E \ar[rr]_{\beta_1} \ar@/^2pc/[rrrr]^{\beta_{r+1}}&& E/\langle p^r Q\rangle \ar[rr]_{\beta_r} && E/\langle Q\rangle}$$
the standard factorization of a cyclic $p^{r+1}$-isogeny into a cyclic $p$-isogeny followed by a cyclic $p^r$-isogeny, where $Q$ is a point of order $p^{r+1}$. 
Let $Q_{r+1}$ and $Q_{r}$ be the unique points such that $e_{\beta_{r+1}}(Q,Q_{r+1})=\zeta_{p^{r+1}}$ and $e_{\beta_{r}}(\beta_1(Q),Q_{r})=\zeta_{p^{r}}$. 
Then $$Q_{r}=pQ_{r+1}.$$
\end{lem}
\proof First of all, we apply the Backing-up Theorem, see \cite[$6.7.11$]{katzmaz}, to the isogeny $\beta_{r+1}$ and to its dual. We have that:
$$\langle \beta_1(Q)\rangle = \ker(\beta_r),\quad \langle p^r Q\rangle = \ker(\beta_1),$$ 
and  
$$\langle \beta^t_r(Q_{r+1})\rangle = \ker(\beta^t_1), \quad \langle p(Q_{r+1})\rangle = \ker(\beta^t_r).$$
For the sake of clarity, we resume all the isogenies we have taken into account in the following diagram:
$$\xymatrix{
\langle Q\rangle\; \ar@{>->}[dr] & \quad && \quad && \quad & \;\; \langle\beta_1(Q)\rangle \ar@{>->}[dl]\\ 
\langle p^r Q\rangle\;\ar@{>->}[r] & E\; \ar@/^/[rrrr]^{\beta_1} \ar@/^/[ddrr]^{\beta_{r+1}} && \quad && E/\langle p^r Q\rangle \ar@/^/[llll]^{\beta_1^{t}}\ar@/^/[ddll]^{\beta_{r}}& \;\langle \beta^t_r(Q_{r+1})\rangle\; \ar@{>->}[l]\\
\quad & \quad && \quad && \quad & \quad\\
\quad & \quad && E/\langle Q\rangle\; \ar@/^/[uull]^{\beta_{r+1}^t}\ar@/^/[uurr]^{\beta_{r}^t} && \quad & \quad\\
\quad & \langle pQ_{r+1}\rangle\; \ar@{>->}[urr] && \quad &&\; \langle Q_{r+1}\rangle \ar@{>->}[ull] & \quad}$$
By definition of the Weil pairing, we have that $\langle Q_{r}\rangle = \ker(\beta^t_r)=\langle pQ_{r+1}\rangle$, so there exists $k\in (\Z/p^r\Z)^\ast$ such that 
$kQ_r = pQ_{r+1}$, since they generate the same cyclic group of order $p^r$. Hence, by bilinearity of the pairing we have 
$$e_{\beta_r}(\beta_1(Q),k Q_{r})=\zeta_{p^{r}}^k= e_{\beta_r}(\beta_1(Q), pQ_{r+1}).$$  
The compatibility of the pairing, see \cite[2.8]{katzmaz} and \cite[III~$8.1$]{silv}, implies that 
$$e_{\beta_r}(\beta_1(Q), pQ_{r+1})=e_{\beta_{r+1}}(Q,pQ_{r+1})=\zeta_{p^{r+1}}^p= \zeta_{p^r}.$$
As $e_{\beta_r}(\beta_1(Q), Q_{r})=\zeta_{p^r}$, then $k=1$. Therefore, $Q_{r}=pQ_{r+1}$.
\endproof

\begin{lem} \label{old2}
Let $n$ and $k$ be positive integers, and let $p$ be a prime such that $n=m p^{r+1}$ with $p$ and $m$ coprime, $r\in\Z_{\geq 1}$ and $m p^{r}>4$. 
Let $\ell$ be a prime not dividing $n$ and such that $1\leq k \leq \ell{+}1$. Let $f \in S(m, p^{r+1}, k)_{\overline{\F}_\ell}$ such that 
$$f(\tate(q), \zeta_{m}, \zeta_{p^{r+1}})= \sum_{j\geq 1} a_j(f) q^{j} \left(\frac{dt}{t}\right)^{\otimes k}\in \overline{\F}_\ell\left[\left[q^p\right]\right]\left(\frac{dt}{t}\right)^{\otimes k},$$
where $\zeta_{m}$ and $\zeta_{p^{r+1}}$ are respectively fixed $m$-th and $p^{r+1}$-th roots of unity. 
Then there exists a unique cuspform $g\in S(m, p^r,k)_{\overline{\F}_\ell}$, such that $B_p^\ast g=f$, where $B_p^\ast$ is the pullback of 
$B_p: X_1(m,p^{r+1})_{\overline{\F}_\ell} \to X_1(m, p^r)_{\overline{\F}_\ell}$ the $p$-th degeneracy map.
\end{lem}
\proof 
Since $m p^{r}>4$ we have the following diagram:
$$\xymatrix{X_1(m,p^{r+1})_{\overline{\F}_\ell} \ar[d]_{B_p} \ar[rr]^*+{w_{\zeta_{p^{r+1}}}} && X_1(m, p^{r+1})_{\overline{\F}_\ell} \ar[d]^{\tilde{\alpha}}\\
X_1(m,p^{r})_{\overline{\F}_\ell} \ar[rr]^*+{w_{\zeta_{p^{r}}}} && X_1(m,p^{r})_{\overline{\F}_\ell}}$$
where $B_p$ is the $p$-th degeneracy map, $w_{\zeta_{p^r}}$ and $w_{\zeta_{p^{r+1}}}$ are Atkin-Lehner involutions and $\tilde{\alpha}:= w_{\zeta_{p^r}}\circ B_p \circ w_{\zeta_{p^{r+1} }}^{-1}$.
Our first goal is to prove that the map $\tilde{\alpha}$ is the forgetful map $\alpha$. Through the moduli interpretation, we can write the previous diagram as follows. 
Let $E$ be an elliptic curve over $S$, an $\overline{\F}_\ell$-scheme, with $P$ and $Q$ respectively points of order $m$ and $p^{r+1}$, we have:
$$\xymatrix @R=2pc @C=1pc {
(E,P,Q) \ar@{|->}[d]_{B_p} \ar@{|->}[rr]^(.4)*+{w_{\zeta_{p^{r+1}}}} && (E/\langle Q\rangle, \beta_{r+1}(P), Q_{r+1}) \ar@{|->}[d]^{\tilde{\alpha}}\\
(E/\langle p^r Q\rangle, \beta_{1}(P), \beta_1(Q)) \ar@{|->}[rr]^*+{w_{\zeta_{p^{r}}}}&& (E/\langle Q\rangle, \beta_r(\beta_1(P)), Q_{r}) }$$
while the forgetful map gives
$$\xymatrix @R=2pc @C=1pc {(E/\langle Q\rangle, \beta_{r+1}(P), Q_{r+1}) \ar@{|->}[rr]^{\alpha} && (E/\langle Q\rangle, \beta_{r+1}(P), pQ_{r+1})},$$
where the maps $\beta_1, \beta_r, \beta_{r+1}$ are isogenies defined by: 
$$\xymatrix{ E \ar[rr]_{\beta_1} \ar@/^2pc/[rrrr]^{\beta_{r+1}}&& E/\langle p^r Q\rangle \ar[rr]_{\beta_r} && E/\langle Q\rangle},$$
so $\beta_r(\beta_1(P))=\beta_{r+1}(P)$. In order to show that $\tilde{\alpha}=\alpha$, it is sufficient to apply Lemma~\ref{isog}: $Q_r = pQ_{r+1}$. 

We have that:
$$w_{\zeta_{p^{r+1}}}(\infty)=w_{\zeta_{p^{r+1}}}(\tate(q), \zeta_{m}, \zeta_{p^{r+1}})= (\tate(q^{p^{r+1}}), \zeta_{m}^{p^{r+1}}, q)=\infty'.$$
There exists a unique Katz cuspform $f' \in S(m, p^{r+1},k)_{\overline{\F}_\ell}$ such that $w_{\zeta_{p^{r+1}}}(f')=f$: by  direct computation at the standard cusp, we have that for all positive integers $j$ $$a_j(f')p^{(r+1)k}=a_j(f),$$ hence 
$f'(\tate(q), \zeta_{m}, \zeta_{p^{r+1}})\in \overline{\F}_\ell\left[\left[q^p\right]\right]\left(dt/t\right)^{\otimes k}$.  

Let us consider the following covering of modular curves:
$$
\xymatrix{
X(\Gamma_1(m),\Gamma(p^{r+1})^{\zeta_{p^{r+1}}-\can})_{\overline{\F}_\ell}
\ar[d]_\gamma^{\left({\begin{smallmatrix}1 & \ast  \\0 & 1\end{smallmatrix}}\right)} 
\ar@/_5pc/[dd]_{\overline{\Gamma_1(p^r)}} && (E,P,Q_{1},Q_{2})\ar@{|->}[d]\\
X(\Gamma_1(m),\Gamma_1(p^{r+1}))_{\overline{\F}_\ell} \ar[d]_{\alpha} && (E,P,Q_{1})\ar@{|->}[d]\\
X(\Gamma_1(m),\Gamma_1(p^{r}))_{\overline{\F}_\ell} && (E,P,p Q_{1})}
$$
where 
$$\overline{\Gamma_1(p^r)}\colon= 
\left\{\left(\begin{matrix}1+p^r\cdot a  & b\\p^r \cdot c & 1+p^r \cdot d \end{matrix}\right)\in \SL_2(\Z/p^{r+1}\Z)\right\}.$$ 

There exists $g \in  S(m,p^{r},k)_{\overline{\F}_\ell}$ such that $B^\ast_p g=f$ if and only if $\gamma^\ast f'$ is $\overline{\Gamma_1(p^r)}$-invariant. Since the set $\left\{\abc{1}{0}{p^r}{1}, \abc{1}{1}{0}{1}\right\}$ 
generates $\overline{\Gamma_1(p^r)}$ by Lemma~\ref{gen}, it is enough to show the invariance for the generators.

The map $\gamma : X(\Gamma_1(m),\Gamma(p^{r+1})^{\zeta_{p^{r+1}} - \can})_{\overline{\F}_\ell} \to X(\Gamma_1(m),\Gamma(p^{r+1}))_{\overline{\F}_\ell}$ corresponds to the action of elements of the form 
$\left({\begin{smallmatrix}1 & \ast  \\0 & 1\end{smallmatrix}}\right)\in \SL_2(\Z/p^{r+1}\Z)$ and $\gamma^\ast f'$ is invariant under the action of such elements by construction. 
The action of $\left({\begin{smallmatrix}1 & 0 \\p^r & 1\end{smallmatrix}}\right)$ fixes 
$\infty'\in X(\Gamma_1(m),\Gamma(p^{r+1})^{\zeta_{p^{r+1}} - \can})_{\overline{\F}_\ell}$ and it acts on $\overline{\F}_\ell\left[\left[q\right]\right]$ as
$$
\left\{
\begin{array}{l}
\mbox{ identity on }\overline{\F}_\ell\\
q \mapsto \zeta_{p^{r+1}}^{p^r} q	=\zeta_p q.
\end{array}
\right.,
$$  
Then, by definition, $\left({\begin{smallmatrix}1 & 0 \\p^r & 1\end{smallmatrix}}\right)$ acts on 
$\overline{\F}_\ell\left[\left[q\right]\right]\left(dt/t\right)^{\otimes k}$ by $t \mapsto t$. 
Since $a_j(f')=0$ for all po\-si\-ti\-ve integers $j$ not divisible by $p$, then $\left({\begin{smallmatrix}1 & 0 \\p^r & 1\end{smallmatrix}}\right)^\ast f'$ has the same $q$-expansion as $f'$ 
and since $X(\Gamma_1(m),\Gamma(p^{r+1})^{\zeta_{p^{r+1}} - \can})_{\overline{\F}_\ell}$ is an irreducible reduced curve, this means that 
$\left({\begin{smallmatrix}1 & 0 \\p^r & 1\end{smallmatrix}}\right)^\ast\gamma^\ast  f'=\gamma^\ast f'$.

The uniqueness of $g\in S(m, p^r,k)_{\overline{\F}_\ell}$ such that 
$B^\ast_p g=f$ follows, as in the previous lemma, by the $q$-expansion principle.\endproof

In the previous lemmas, we did not consider the action of diamond operators. First of all, let us recall briefly how this operators are defined.																																										
For $n > 4$, and for all $d \in (\Z/n\Z)^\ast$, we define an automorphism
\begin{eqnarray*}
r_d:X(n)_{\overline{\F}_\ell} & \to & X(n)_{\overline{\F}_\ell}\\\nonumber
(E, P) &\mapsto & r_d(E, P) = (E, dP),\nonumber
\end{eqnarray*} 
for all generalized elliptic curves $E$ together with a $\Gamma_1(n)$-structure $P$. 
For every $d \in (\Z/n\Z)^\ast$, the diamond operator $\langle d \rangle$ on $M(\Gamma_1(n),k)_{\overline{\F}_\ell}$  is defined as the automorphism of $M(\Gamma_1(n),k)_{\overline{\F}_\ell}$ 
induced by pullback via the automorphism $r_d$. 
If we are dealing with Katz cuspforms with a character $\epsilon:(\Z/n\Z)^\ast \to \overline{\F}_\ell^\ast$, we are imposing that for all $d \in (\Z/n\Z)^\ast$, the diamond operator $\langle d \rangle$ acts as $\epsilon(d)$. 

\begin{lem}
Let $n$ and $k$ be positive integers, let $p$ be a prime dividing $n$ such that $n/p>4$. Let $\ell$ be a prime not dividing $n$ such that $1\leq k \leq \ell{+}1$. 
Let $\epsilon:(\Z/n\Z)^\ast \to \overline{\F}_\ell^\ast$ and $\chi: (\Z/(n/p)\Z)^\ast\to \overline{\F}_\ell^\ast$ be characters. 
Let $f\in S(n, k, \epsilon)_{\overline{\F}_\ell}$ and $g\in S(n/p, k, \chi)_{\overline{\F}_\ell}$ such that $B^\ast_p (g)=f$. Then
$$
\chi= \res_{\Z/n\Z}^{\Z/(n/p)\Z}(\epsilon).
$$
\end{lem} 
\proof
Let $E/S$ be an el\-liptic curve over an $\overline{\F}_\ell$-scheme $S$, let $P$ be a point of order $n$, and let $B_p$ be the degeneracy map $B_p: X_1(n)_{\overline{\F}_\ell} \to X_1(n/p)_{\overline{\F}_\ell}$ 
defined, as in (\ref{dege}), through the isogeny $\beta$ then 
$$(B_p^\ast g)(E,P)=\beta^\ast (g(E/\left\langle (n/p) P\right\rangle, \beta(P))).$$
The action of the diamond operator $\left\langle d\right\rangle$ for $d\in \Z/(n/p)\Z$ is given by
\begin{align*}
\left\langle d\right\rangle(B_p^\ast g)(E,P)&=(B_p^\ast g)(E,dP)=\beta^\ast (g(E/\left\langle (n/p) P\right\rangle, \beta(dP)))=\\
&=\beta^\ast (g(E/\left\langle (n/p) P\right\rangle, d\beta(P)))=\\
&=\beta^\ast (\chi(d) g(E/\left\langle (n/p) P\right\rangle, \beta(P)))= \chi(d) (B_p^\ast g)(E,P);
\end{align*}
hence, for all $d\in \Z/(n/p)\Z$ we have $\left\langle d\right\rangle(B_p^\ast g)=\chi(d) (B_p^\ast g)$. Since $B^\ast_p (g)=f$, then for all $d\in \Z/(n/p)\Z$ it follows that 
$\left\langle d\right\rangle f =\chi(d)f =\epsilon(d) f$. Therefore, for all $d\in \Z/(n/p)\Z$ the equality $\chi(d)=\epsilon(d)$ holds.
\endproof

If we take this into account, we have a generalization of Lemma~\ref{old1} and Lemma~\ref{old2} for forms in $S(n, k, \epsilon)_{\overline{\F}_\ell}$: 
\begin{lem}\label{cold2}
Let $n$ and $k$ be positive integers. Let $p$ be a prime such that $n$ factors as $m p^{r+1}$ with $p$ and $m$ coprime, $r\in\Z_{\geq 0}$ and $m p^{r}>4$. 
Let $\ell$ be a prime not dividing $n$ and such that $1\leq k \leq \ell{+}1$. 
Let $\epsilon:(\Z/n\Z)^\ast \to \overline{\F}_\ell^\ast$ be a character. 
Let $f\in S(m, p^{r+1}, k, \epsilon)_{\overline{\F}_\ell}$ such that 
$$f(\tate(q), \zeta_{m}, \zeta_{p^{r+1}})= \sum_{j\geq 1} a_j(f) q^{j} \left(\frac{dt}{t}\right)^{\otimes k}\in \overline{\F}_\ell\left[\left[q^p\right]\right]\left(\frac{dt}{t}\right)^{\otimes k},$$
where $\zeta_{m}$ and $\zeta_{p^{r+1}}$ are respectively fixed roots of unity of order $m$ and $p^{r+1}$. 
Then there exists a unique form $g\in S(m, p^r,k, \epsilon')_{\overline{\F}_\ell}$, such that $B_p^\ast g=f$, where $\epsilon'= \res^{\Z/(mp^r)\Z}_{\Z/n\Z}(\epsilon)$.
\end{lem}

In Lemma~\ref{old1} and~\ref{old2}, as well as in the corollaries above, we assume that the prime $p$, dividing the level $n$, is such that the quotient $n/p >4$. 
This hypothesis implies that all the objects used in the proofs are modular curves.  Without such hypothesis, we should use the theory of algebraic stacks, see \cite[Theorem~$1.2.1$]{con}.
Anyway, we can apply level raising, introducing a new prime in the level so that the hypotheses of Lemma~\ref{cold2} are satisfied. 
Hence, after applying Lemma~\ref{cold2}, we can observe that the form obtained comes from the starting level, it belongs to the span of the degeneracy maps at the prime introduced in the level. 
Therefore the following  holds:
\begin{cor}\label{cold}
Lemma~\ref{cold2} holds with no assumption on $n/p$.
\end{cor}
Moreover, we can also generalize the previous result, iterating through divisors and equating $q$-expansions:
\begin{lem}
Let $n$ and $k$ be positive integers. Let $d$ be a positive integer such that $n$ factors as $m d$ with $d$ and $m$ coprime. 
Let $\ell$ be a prime not dividing $n$ and such that $1\leq k \leq \ell{+}1$. 
Let $\epsilon:(\Z/n\Z)^\ast \to \overline{\F}_\ell^\ast$ be a character. 
Let $f\in S(m, d, k, \epsilon)_{\overline{\F}_\ell}$ such that 
$$f(\tate(q), \zeta_{m}, \zeta_{d})= \sum_{j\geq 1} a_j(f) q^{j} \left(\frac{dt}{t}\right)^{\otimes k}\in \overline{\F}_\ell\left[\left[q^d\right]\right]\left(\frac{dt}{t}\right)^{\otimes k},$$
where $\tate(q)$ is the Tate curve over $\overline{\F}_\ell((q))$, $\zeta_{m}$ and $\zeta_d$ are fixed roots of unity of order $m$ and $d$. 
Then there exists a unique form $g\in S(m,k, \epsilon')_{\overline{\F}_\ell}$, such that $B_p^\ast g=f$, where $\epsilon'= \res(\epsilon)$.
\end{lem}
The proof of this statement consists in the iteration of applications of Lemma~\ref{cold2} and Corollary~\ref{cold}, after factorizing $d$.
\section{Proof of Theorem~\ref{rid}}

%
%
%
%
%

\proof 
The operator $T_r$ commutes with all $T_p$, and it acts on $V$. Let $v$ be the dimension of $V$. 
By hypothesis $v>1$, hence $V$ contains a non-zero eigenvector $f$ for all $T_p$. In particular $a_1(f)\neq 0$. 

The subspace $V_1\subset V$ consisting of the $g$ in $V$ with $a_1(g)=0$ is of dimension $v-1>0$. 
Every element $g\in V_1$ is such that $0=a_1(g)=a_1(T_j(g))$ for all integer $j$ not divisible by $m$. 
This means that $g$ is an eigenform with $q$-expansion in $\F\left[\left[q^r\right]\right]\left((dt)/t\right)^{\otimes k}$ at the standard cusp. 

Let us first suppose that $r\neq \ell$. Then we are in the hypotheses of Lemma~\ref{cold2} and Corollary~\ref{cold}, so there exists a non zero mod~$\ell$ 
cuspform $g' \in S(n/r, k, \epsilon')_{\overline{\F}_\ell}$ such that $B^\ast_r g'=g$, where $\epsilon'=\res(\epsilon)$ and $B_r^\ast$ is the pullback of the degeneracy map $B_r$. 
This means that $V_1= B^\ast_r V_2$ where $V_2$ is the eigenspace of 
$S(n,k, \overline{\epsilon})_{\overline{\F}_\ell}$ associated to the system of eigenvalues $\{(f(T_p), f(\left\langle p\right\rangle)) \}$ for $p$ prime, $p$ different from $r$.
Hence the system of eigenvalues occurs in level $n/r$.

Let us now suppose that $r$ is equal to $\ell$. As proved in \cite[Proposition~6.2]{edix1}, the only two following cases occur. First case: $k=\ell$. The form $g$, constructed as before, is in the image of the Frobenius $F$, i.e.\ there exists $g'\in S(n, 1, \epsilon)_\F$ such that $g=(g')^\ell=F(g')$. 
This means that the system of eigenvalues occurs already in weight $1$. Second case: $k=\ell +1$. We have that $g=A_\ell g'$ where $A_\ell$ is the Hasse invariant of weight $\ell-1$ and  $g'\in S(n, 2, \epsilon)_\F$, hence the system of eigenvalues occurs in weight $2$.
\endproof

From the previous theorem, recursively applying the previous argument to the subspace $V_1$ constructed in the proof above, we deduce the following:

\begin{cor}
In the same hypotheses of Theorem~\ref{rid}, if $r\neq \ell$ then the dimension of $V$ is at most $a+1$ where $a$ is such that $r^a \nc(\rho_f)=n$.
\end{cor}

\begin{ex}
In \cite{diso}, Dieulefait and Soto gave an example of a mod~$5$ eigenvalue system attached to a modular form $f$ of level $\Gamma_0(1406)$ and weight $2$ which satisfies the hypotheses of Theorem~\ref{rid}: 
there exist a form $g$ at the same level and weight $2$ whose $q$-expansion coincide with the one of $f$ away from $19$, we have indeed $1=f(T_{19})=-g(T_{19})$.  
It is easy to check that the mod~$5$ eigenvalue system arises from a modular form of level $\Gamma_0(74)$ and weight $2$: there are two Galois orbits of eigenforms in level $\Gamma_0(74)$ and weight $2$, and the eigenvalue system arises by reduction of a characteristic zero eigenvalue system in the orbit with Hecke eigenvalue field is $\Q(\sqrt{5})$.
\end{ex}

\section*{Acknowledgements}
I would like to thank Bas Edixhoven for the long discussions about the material in this article. 
I also am very grateful to Gabor Wiese for all the fruitful conversations about this topic.  

\bibliographystyle{annotate}
\bibliography{biblio}

\end{document}